\documentclass[letterpaper,11pt,openany,oneside]{article}
\usepackage[T1]{fontenc}
\usepackage[latin1]{inputenc}
\usepackage{epsfig}
\usepackage{amsmath,amssymb}
\usepackage[active]{srcltx}
\usepackage{dsfont}
\usepackage[titles]{tocloft}
\usepackage{amssymb}
\usepackage{tikz}
\usetikzlibrary{decorations.pathreplacing}
\numberwithin{equation}{section}
\newtheorem{theoreme}{Theorem}[section]

\newtheorem{remarque}[theoreme]{Remark}

\newtheorem{lemme}[theoreme]{Lemma}

\newtheorem{definition}[theoreme]{Definition}

\newtheorem{proposition}[theoreme]{Proposition}

\newcommand{\RR}{\ensuremath{\mathbb R}}

\newcommand{\NN}{\ensuremath{\mathbb N}}

\newcommand{\De}{\Delta}

\newcommand{\de}{\delta}

\title{Uniform value for some nonexpansive optimal control problems with general evaluations}
\author{Xiaoxi Li\thanks{Economics and Management School, Wuhan University, Luojia Hill, Wuchang, 430072 Wuhan, China. Email: xxleewhu@gmail.com.}}

\date{\today }

\begin{document}
\maketitle
\begin{abstract}  
We study the long-run properties of optimal control problems in continuous time, where the running cost of a control problem is evaluated by a probability measure over $\RR_+$. \cite{Li_2016} introduced an asymptotic regularity condition for a sequence of probability measures to study the limit properties of the value functions with respect to the evaluation. In the particular case of $t$-horizon Ces\`aro mean or $\rho$-discounted Abel mean, this condition implies that the horizon $t$ tends to infinity or the discount factor $\rho$ tends to zero. For the control system defined on a compact domain and satisfying some nonexpansive condition, \cite{Li_2016} proved the existence of general limit value, i.e. the value function uniform converges as the evaluation becomes more and more regular. Within the same context, we prove the existence of general uniform value, i.e. for any $\varepsilon>0$, there is an optimal control that guarantees the general limit value up to $\varepsilon$ for all control problems where the cost is evaluated by sufficiently regular probability measures.   
\end{abstract}



\section{Introduction}
Let $U$ be a compact subset of a separable metric space. A \textit{control} ${\bf{u}}$ is a measurable function from $\RR_+$ to $U$. Denote by $\mathcal{U}$ the set of all controls. We consider the following control system:
\begin{equation} \label{eq:dynamics}
{\bf y}'(t)=f\big({\bf y}(t), {\bf u}(t)\big), \ {\bf y}(0)=y_0.
\end{equation}
where $f:\RR^d\times U \to \RR^d$, and $y_0\in\RR^d$ is the initial state. Let $g:\RR^d\times U\to \RR$ be the running cost function. We make the following assumptions on $f$ and $g$ throughout the article: 
\begin{equation}\label{ex:assumptions}
\begin{cases}
 \text{the function}  \  g:\RR^d\times U\to\RR \text{ is Borel measurable and bounded}; \\
 \text{the function} \ f:\RR^d\times U\to\RR^d \text{ is Borel measurable, and satisfies:}\\
(*).\ \exists L\geq 0, \forall (y,\overline{y})\in\RR^{2d}, \forall u\in U, ||f(y,u)-f(\overline{y},u)||\leq L||y-\overline{y}||,\\
(**).\ \exists a>0, \forall (y,u)\in\RR^d\times U, ||f(y,u)||\leq a(1+||y||).
\end{cases}
\end{equation}
Then for a given initial point $y_0\in \RR^d$, any control ${\bf{u}}\in\mathcal{U}$ defines a unique absolutely continuous solution to Equation (\ref{eq:dynamics})  on $\RR_+$, which we denote by ${\rm \textbf{y}}(t,\bf{u},y_0)$.   

Denote by $\mathcal{J}=\langle U, g,f\rangle$ the optimal control problem described above. $\De(\RR_+)$ is the set of probability measures on $\RR_+$ and any $\theta\in\De(\RR_+)$ is called  an \textit{$evaluation$}. The $\theta$-value of the control problem $\mathcal{J}$ is defined as:
\begin{equation}\label{eq:V-theta}
V_\theta(y_0)=\inf_{{\bf{u}}\in\mathcal{U}}\gamma_\theta(y_0,{\bf{u}}), \text{ \ with \ } \gamma_\theta(y_0,{\bf{u}})=\int_{0}^{+\infty}g\big({\rm \textbf{y}}(t,{\bf{u}},y_0),{\bf{u}}(t)\big) d\theta(t). 
\end{equation}

Specific evaluations and the corresponding value functions include

\noindent \textit{Ces\`aro  mean}: \ $\forall t>0$, $\overline{\theta}_t$ has a density $s\mapsto f_{\overline{\theta}_t}(s)=\frac{1}{t}1_{[0,t]}(s)$, and the $t$-horizon value is
\begin{equation*}
V_{\overline{\theta}_t}(y_0)=\inf_{{\bf{u}}\in \mathcal{U}} \frac{1}{t}\int_{0}^{t} g\big({\bf{y}}(s,{\bf{u}},y_0),{\bf{u}}(s)\big) d s
\end{equation*}
\noindent \textit{Abel mean}: \ $\forall \rho\in(0,1]$, $\theta_\rho$ has a density $s\mapsto f_{\theta_\rho}(s)=\rho e^{-\rho s}$, and the $\rho$-discounted value is
\begin{equation*}
V_{\theta_\rho}(y_0)=\inf_{{\bf{u}}\in \mathcal{U}}\int_{0}^{+\infty} \rho e^{-\rho s}g\big({\bf{y}}(s,{\bf{u}},y_0),{\bf{u}}(s)\big) d s
\end{equation*}

In this article, we are interested in the long-run properties of the control problem $\mathcal{J}$. In particular cases, we look at the convergence of $V_{\bar{\theta}_t}$ as $t$ tends to infinity, the convergence of $V_{\theta_\rho}$ as $\rho$ tends to zero, and their equality (the \textit{asymptotic approach}). In case of convergence, the \textit{uniform approach} further asks for each $\varepsilon>0$ the existence of  $\varepsilon$-optimal control for all control problems with sufficiently large $t$ and/or with sufficiently small $\rho$. The uniform approach emphasizes the robustness of the optimal control with respect to a long but unknown duration of the control problem.  

In the literature, the asymptotic approach has been extensively exploited. Most of the results are obtained by assuming certain ergodic condition or controllability condition (cf. \cite{Alvarez_2010, Arisawa_97, Arisawa_98, Bensoussan_88, Gaitsgory_86}). A different approach making no ergodic/controllability assumption is initialized by \cite{Quincampoix_2011}. Notably, the obtained limit value function (whenever exists) is in general dependent of the initial state. To be more precise, they pose a \textit{compact nonexpansive} assumption on the control problems, i.e. the control system is defined on a compact domain and satisfies certain mild nonexpansive condition. It is proven that for any compact nonexpansive control problem, $V_{\bar{\theta}_t}$ converges uniformly (in the initial state) as $t$ tends to infinity, and moreover for any $\varepsilon>0$, $\varepsilon$-optimal control exists for all $t$-horizon problems with sufficiently large $t$. 

This nonexpansive approach for the asymptotic study of the average value in optimal control has recently been advanced in the literature toward different directions: \cite{Cannarsa_2015} for optimal control and differential games linked to PDE techniques,
\cite{Buckdahn_2014} and \cite{Goreac_2015b} for stochastic control system, \cite{Quincampoix_2015} for singularly perturbed control system.  

\bigskip

On the other hand, several Tauberian-type results are established in various contexts of optimal control (\cite{Oliu_2013} for deterministic optimal control, \cite{Buckdahn_2014} for stochastic optimal control, \cite{Khlopin_2015} for differential games), which state that the uniform convergence of $t$-horizon value as $t$ tends to infinity is equivalent to that of $\rho$-discounted value as $\rho$ tends to zero, and both limits are equal in case of convergence.   

Motivated by those Tauberian-type results, the asymptotic study of the optimal control problems with general evaluations aims at defining a unique limit value function that is independent of the particular chosen evaluation sequence (for example, to be Ces\`aro mean or Abel mean). For this aim, \cite{Li_2016} introduced the notion of \textit{$s$-total variation} to define an asymptotic regularity condition for evaluation sequences. This condition is used to study the asymptotic behavior of $V_\theta$ as $\theta$ becomes more and more regular, and in particular, it submerges the asymptotic conditions "$t\to+\infty$" and "$\rho\to 0$" as particular cases. In the context of compact nonexpansive control problems, \cite{Li_2016} proved the existence of \textit{general limit value} (cf. Definition \ref{def:glv}), i.e. there is uniform convergence of $V_{\theta}$ as $\theta$ becomes more and more regular. 

\cite{Goreac_2015b} also exploited the asymptotic appoach for control problems with a larger family of evaluations. In his paper, only a family of probability measures that are absolutely continuous with respect to the Lebesgue measure is considered, and a different asymptotic regularity condition from that in \cite{Li_2016}  is introduced for the study.  

\bigskip

Compared to the large body of literature on the asymptotic approach for control problems, results via uniform approach are relatively rare apart from \cite{Quincampoix_2011}. Our article follows the research lines of \cite{Quincampoix_2011} and \cite{Li_2016} to conduct a uniform analysis for compact nonexpansive control problems with general evaluations. Using the asymptotic regularity condition in \cite{Li_2016}, we prove the existence of \textit{general uniform value} (cf. Definition \ref{def:glv-guv}), i.e. for any $\varepsilon>0$, there is $\varepsilon$-optimal control to guarantee the general limit value for all $\theta$-evaluated problems with $\theta$ sufficiently regular. This generalizes the result of uniform approach in \cite{Quincampoix_2011} from $t$-horizon average to general evaluations. 

The proof for the existence of the robust optimal control is constructive.  In \cite{Quincampoix_2011}, to obtain an $\varepsilon$-optimal control for all long but unknown duration, the construction is by blocks such that on each of them, the average payoff is above the limit value, and the value at the starting state is nondecreasing. New difficulty arises in our context since now the optimal control aims at doing well at the same time against a much greater family of evluations, and one can not expect the construction in \cite{Quincampoix_2011} to work as well. To achieve this and especially to obtain a convexity condition needed for the proof of a minmax theorem \big(cf. Equation (\ref{eq:minmax})\big), we introduce randomization on the control.  

Meanwhile, there are some optimality criteria in the infinite-horizon control literature that are related to the uniform approach, such as the \textit{undiscounted criterion} \big(cf. Definition \ref{def:undiscounted-v}\big), where a total cost is defined as the limiting average cost $\limsup_{t\to\infty}\frac{1}{t}\int_{0}^t g\big({\rm \textbf{y}}(s,{\bf{u}},y_0),{\bf{u}}(s)\big)d s $,  or the \textit{weighted average criterion} (cf. Definition \ref{def:weighted-average}), where a total cost is defined as a weighted average of the limiting average cost and the $\rho$-discounted cost. We shall demonstrate that compared with these optimal criteria for long-run control problems, our uniform approach with general evaluations  is rather strong (cf. Proposition \ref{prop:hat-V} and Proposition \ref{prop:weighted-V}).

\bigskip

The organization of the paper is as follows. In Section \ref{sec:prelim}, we present the basic model and introduce the main poblems under study. The main result is presented in Section \ref{sec:mainresult}. We first illustrate it by an example and then discuss its implications for other long-run optimality criteria. The proof of the main result is given in Section \ref{sec:proof}, which consists of two parts: some preliminary results concerning random controls and the main part of the proof.

\section{Preliminaries}\label{sec:prelim}

 Let $\mathcal{J}=\langle U, g,f\rangle$ be optimal control problem defined by Equation (\ref{eq:dynamics}) for which Condition (\ref{ex:assumptions}) is assumed. Let $\theta\in\De(\RR_+)$ be an evaluation for the control problem $\mathcal{J}$, and let $V_\theta(y_0)$ be the $\theta$-value of $\mathcal{J}$ defined by Equation (\ref{eq:V-theta}). The following notion is introduced in \cite{Li_2016} for the aim of defining a regularity for $\theta$.   

\begin{definition} \label{def:LTC} For any $s\geq 0$, the \textbf{$s$-total variation} of an evaluation $\theta$ is: $$TV_s(\theta)=\max_{Q\in\mathcal{B}(\RR_+)}|\theta(Q)-\theta(Q+s)|,$$
where $\mathcal{B}(\RR_+)$ is the family of Borel subsets in $\RR_+$.   
\end{definition}
 
Indeed, we use $\sup_{0\leq s\leq S} TV_s(\theta)$ for some $S>0$ to measure the regularity of any evaluation $\theta$, and we shall conduct the asymptotic analysis and the uniform analysis for the control problem $\mathcal{J}$ as $\sup_{0\leq s\leq S} TV_s(\theta)$ becomes sufficiently small. The reader is refered to Section 3.3 in \cite{Li_2016} for an extensive discussion on this regularity condition. In particular, the vanishing  condition "$\sup_{0\leq s\leq S} TV_s(\theta)\to 0$ for some $S>0$" corresponds to 
either "$t\to +\infty$" for the Ces\`aro mean or "$\rho\to 0$" for the Abel mean.

\bigskip 

The following formal definition is concerned with the asymptotic analysis of control problems with general evaluations.  

\begin{definition}\label{def:glv}
The optimal control problem $\mathcal{J}=\langle U,g,f\rangle$ has a \textbf{general limit value} given by some function $V$ defined on $\RR^d$ if: for each $\varepsilon>0$ there is some $\eta>0$ and $S>0$ such that:
$$\forall \theta\in\De(\RR_+), \ \  \Big(\sup_{0\leq s\leq S}TV_s(\theta)\leq \eta \Longrightarrow \big(\forall y_0\in \RR^d, \ |V_\theta(y_0)-V(y_0)|\leq \varepsilon\big)\Big).$$
\end{definition}

We restrict ourselves to the study of the following class of control problems.
 
\begin{definition} \label{def:compact-nonexpan-control} The optimal control problem $\mathcal{J}=\langle U,f,g\rangle$ is \textbf{compact nonexpansive} if: \\ 
(A.1) the control dynamic (\ref{eq:dynamics}) has a compact invariant set $Y\subseteq\RR^d$, i.e. ${\bf{y}}(t,{\bf{u}},y_0)\in Y$, $\forall t\geq 0$ for all  ${\bf{u}}\in\mathcal{U}$ and $y_0\in Y$;\\ 
(A.2) the running cost $g$ does not depend on $u$ and is continuous in $y\in Y$;\\
(A.3) the control dynamic (\ref{eq:dynamics}) is nonexpansive, i.e., 
\begin{align}\label{eq:non-expansive} \forall y_1,y_2\in\RR^d,\ \sup_{a\in U}\inf_{b\in U}\big\langle y_1-y_2,f(y_1,a)-f(y_2,b)\big\rangle\leq 0.
\end{align}
\end{definition}

\begin{remarque} The nonexpansive condition (A.3) implies the following property (cf.  Proposition 3.7 in \cite{Quincampoix_2011}): $\forall y_1, y_2\in \RR^d$, $\forall {\bf u} \in \mathcal{U}$, $\exists{\bf v}\in\mathcal{U}$ s.t.
$$|| {\bf y}(t, {\bf u}, y_1)-  {\bf y}(t, {\bf v}, y_2)||\leq ||y_1-y_2||,\  \forall t\geq 0.$$  
\end{remarque}

Consider a control problem $\mathcal{J}=\langle U,g,f\rangle$ compact nonexpansive with the invariant set $Y$. \cite{Li_2016} (cf. Corollary 4.5) proved that $\mathcal{J}$ has a general limit value characterized by the following function:
\begin{align} \label{eq:V*}
V^*(y_0)=\sup_{\theta\in\De({\RR_+})}\inf_{s\in\RR_+} \inf_{{\bf{u}}\in\cal{U}}\int_{0}^\infty g\Big({\rm \textbf{y}}(t+s, {\bf{u}},y_0)\Big)d \theta(t), \ \  \forall y_0\in Y.
\end{align}
 
\begin{proposition}\label{prop:LQR}  \cite{Li_2016} \  Let $\mathcal{J}=\langle U,f,g\rangle$ be compact nonexpansive with invariant set $Y$. Then the general limit value of  $\mathcal{J}$ exists and is equal to $V^*$. 
\end{proposition}
 
\begin{remarque}
As opposed to most papers in the control literature, say for example \cite{Alvarez_2010, Arisawa_97, Arisawa_98, Bensoussan_88, Gaitsgory_86},  no ergodic condition or controllability condition is assumed for compact nonexpansive control problems. This can be underlined by the fact that in general the limit function $V^*$ defined in (\ref{eq:V*}) depends on the initial state $y_0$.
\end{remarque}

\bigskip 

We study in this paper a value notion that is stronger than the general limit value, namely the \textit{general uniform value} (the "uniform approach"), which requests the existence of robust controls that is approximately optimal for all $\theta$-evaluated conctrol problems with $\theta$ sufficiently regular (i.e., $\sup_{0\leq s\leq S} TV_s(\theta)$ sufficiently small for some fixed $S>0$).

For the aim of a formal analysis of the uniform approach, we first introduce the notion of \textit{random control}. 

\begin{definition} A \textbf{random control} is a pair $\big((\Omega,\mathcal{B}(\Omega),\lambda), {\bf{u}}\big)$, where $(\Omega,\mathcal{B}(\Omega), \lambda)$ is some standard Borel probability space and ${\bf{u}}:Y\times\Omega\times[0,+\infty)\to U$ is a Borel measurable mapping (we have extended the definition of a random control to dependent on the initial state for later use).
\end{definition}

Denote by $\tilde{\mathcal{U}}$ the set of all random controls, which is convex in the following sense. For any ${\bf u}_1, {\bf u}_2\in \tilde{\mathcal{U}}$ and $\alpha\in[0,1]$, define ${\bf u}=\alpha {\bf u}_1 + (1-\alpha)  {\bf u}_2$  to take the value of ${\bf u}_1$ with probability $\alpha$ and of ${\bf u}_2$ with probability $(1-\alpha)$. It is easy to construct a product probability space  such that ${\bf u}$ is in $\tilde{\mathcal{U}}$ .

We might shortly write $\Omega$ for the triple $(\Omega,\mathcal{B}(\Omega), \lambda)$. Let $(\Omega,{\bf{u}})$ be a random control, then for any initial point $y_0\in Y$ and any $\omega\in \Omega$, ${\bf{u}}_\omega(y_0, \cdot):={\bf{u}}(y_0, \omega, \cdot)$ defined from $[0,+\infty)$ to $U$ is a (pure) control in $\mathcal{U}$, which we denote by ${\bf{u}}_\omega(y_0)$. The Borel probability space $\Omega$ serves as a random device for the controller to choose a pure control in $\mathcal{U}$.

The expected $\theta$-evaluated payoff induced by any random control $(\Omega, {\bf{u}})$ and initial point $y_0$ is denoted by $$\gamma_\theta(y_0,{\bf{u}})=\int_{\Omega}\gamma_\theta\big(y_0, {\bf{u}}_\omega(y_0)\big)d \lambda(\omega)=\int_{\Omega}\int_{[0,+\infty)} g\Big({\rm \textbf{y}}(t,{\bf{u}}_{\omega}(y_0),y_0)\Big)d \theta(t)d \lambda(\omega),$$
and the expected $\theta$-value in random controls is $\tilde{V}_{\theta}(y_0)=\inf_{{\bf{u}}\in\mathcal{\tilde{\mathcal{U}}}}\gamma_\theta(y_0, {\bf{u}})$. The payoff function $\gamma_\theta(y_0,\cdot)$ is affine \footnote{To see this point, consider for example the Borel probability space to be $\big([0,1],\mathcal{B}([0,1]),\lambda \big)$ where $\lambda$ is the Lebesgue measure. We take any ${\bf{u}}^1, {\bf{u}}^2:Y\times [0,1]\times[0,+\infty)\to U$ two random controls, any $\alpha\in[0,1]$, and let ${\bf{u}}^3:Y\times [0,1]\times[0,+\infty)\to U$ be the random control as a convex combination of ${\bf{u}}^1$ and ${\bf{u}}^2$ with coefficient $\alpha$: for any $y_0\in Y$, ${\bf{u}}^3(y_0,\omega, t)={\bf{u}}^1(y_0, \frac{\omega}{\alpha}, t)$ for $\omega\in[0,\alpha]$ and ${\bf{u}}^3(y_0,\omega, t)={\bf{u}}^2(y_0,\frac{\omega-\alpha}{1-\alpha}, t)$ for $\omega\in(\alpha,1]$. Using the change of variables, we obtain $\gamma_{\theta}(y_0, {\bf{u}}^3)=\alpha\cdot\gamma_\theta(y_0, {\bf{u}}^1)+ (1-\alpha)\cdot\gamma_\theta(y_0, {\bf{u}}^2)$.} in ${\bf{u}}$, thus the value function in random controls is the same as that in pure controls, that is, $\tilde{V}_{\theta}(y_0)=V_{\theta}(y_0)$ for all $y_0\in Y$ and $\theta\in\De(\RR_+)$. 

\begin{definition} \label{def:glv-guv} The optimal control problem $\mathcal{J}=\langle U,f,g\rangle$ has a \textbf{general uniform value} if the general limit value exists, say equal to $V^*$, and moreover for each $\varepsilon>0$ there is some $\eta>0$, $S>0$ and a random control ${\bf{u}}\in\tilde{\mathcal{U}}$ such that:
$$\forall \theta\in\De(\RR_+), \ \  \Big(\sup_{0\leq s\leq S}TV_s(\theta)\leq \eta \Longrightarrow \big(\forall y_0\in Y, \ \gamma_\theta(y_0,{\bf{u}})\leq V^*(y_0)+\varepsilon\big)\Big).$$ 
\end{definition}

The control ${\bf{u}}$ appearing in the above definition is called an \textit{$\varepsilon$-optimal control} for the control problem $\mathcal{J}$.

\section{Main result}\label{sec:mainresult}

Our main result of this article is the following

\begin{theoreme}\label{thm:guv} \ Assume that the optimal control problem $\mathcal{J}=\langle U,f,g\rangle$ is compact nonexpansive. Then it has a general uniform value $V^*$.
\end{theoreme}

\begin{remarque} \cite{Quincampoix_2011} proved the existence of pure $\varepsilon$-optimal control for the compact nonexpansive control problems when the running cost is evaluated by Ces\`aro means.  Our result generalizes it to general evaluations. Nevertheless, new difficulty arises and our proof techniques differ from \cite{Quincampoix_2011} and is partially inspired by \cite{Renault_2013}, where some analogous results in a discrete-time framework  are  established.

It is unknown whether pure $\varepsilon$-optimal control exists for compact nonexpansive control problems with general evaluations. 
\end{remarque}

The rest of this section is devoted to comments of our main result. We first present in Subsection \ref{subsec:example} a toy example of control problem that is compact nonexpansive thus general uniform value exists. Next in Subsection \ref{subsec:implications}, we link the general uniform value to other long-run optimal criteria in the control literature.  

\subsection{Illustration of Theorem \ref{thm:guv} by a toy example}\label{subsec:example} This subsection contains a simple example  of economic development. We show that its associated control problem is compact nonexpansive, thus according to Theorem \ref{thm:guv} the general uniform value exists.

A society is concerned with its long-run economic development -- which is denoted by a state variable $x_1$, and also with the accompanied pollution impact -- which is denoted by a state variable $x_2$. Let $\pi(x_1)$ be the instantaneous benefit (assumed to be continuous) of the development and let $c(x_2)$ be the instantaneous cost (assumed to be continuous) of the pollution. $u\in[0,1]=U$ is the control variable, which can be interpreted as the policy of develeopment, such as the number of vehicles allowed to enter into the metropolis per day, or the volume of carbon released into the atmosphere per year, etc. We denote $X=(x_1, x_2)$, and the dynamic of $X(t)$ is governed by ($A, B>0$ constant)
\[X'(t)=f\big(X(t), u(t)\big)= \left( \begin{array}{c}
u(t)  \big(A-x_1(t)\big)\\
u(t)\big(B-x_2(t)\big)  \end{array} \right), \ \ \forall t\geq 0.\] 
Starting from $X(0)=(0,0)$, we have $x'_1(t)\geq 0, x'_2(t)\geq 0$ and $x_1(t)\leq A,  x_2(t)\leq B$ for all $t\geq 0$. $u(t)$ is the control level at time $t$, which promotes a positive growth of both the develepment $x_1(t)$ and the pollution $x_2(t)$. 

Let us define $g(X)=c(x_2)-\pi(x_1)$ to be the runnning cost function, then the associated optimal control problem $\mathcal{J}=\langle U,f,g\rangle$ satisfies the assumptions in (\ref{ex:assumptions}). Moreover, we verify that $\mathcal{J}$ is compact nonexpansive thus it has a general uniform value:

\begin{itemize}
\item $[0,A]\times [0, B]$ is its compact invariant set; 
\item $g(X,u)=g(X)$ does not depend on $u$ and is continuous in $X$;
\item the control dynamic $X'(t)=f\big(X(t), u(t)\big)$ is nonexpansive: let $X=(x_1, x_2)$ and $Z=(z_1, z_2)$, for any $u\in [0,1]$, take $v=u$. Then we obtain:
\[ \Big\langle X-Z, f(X, u)-f(Z, u)\Big\rangle= \Bigg\langle \left( \begin{array}{c}
x_1-z_1\\
x_2-z_2  \end{array} \right), u\left( \begin{array}{c}
 \big(-x_1+z_1\big)\\
\big(-x_2+z_2\big)  \end{array} \right) \Bigg\rangle\leq 0.\] 
\end{itemize}

\subsection{Implications of the general uniform value to other optimality criteria}\label{subsec:implications}

In this subsection, we link the general uniform value to other solution criteria such as the \textit{undiscounted criterion} (limiting avarage payoff, cf. \cite{Carlson_91}) and the \textit{weighted average criterion} (cf.  \cite{Filar_92, Haurie_2002}). As the asymptotic/uniform approach, both the undiscounted value and the limit value of the weighted average value (when the discount factor tends to zero) are also optimality criteria for long-run control problems. We show that the general uniform value is actually a stronger value notion.

\bigskip 

Below we assume that the control problem $\mathcal{J}=\langle U,f,g\rangle$ under consideration is compact nonexpansive, and let $Y\subseteq \RR^d$ be its invariant set. Following Theorem \ref{thm:guv}, $\mathcal{J}$ has a general uniform value $V^*$.
 
\begin{definition}\label{def:undiscounted-v}  The \textbf{undiscounted value} $\hat{V}$ of the control problem $\mathcal{J}$ with limiting average payoff is defined as:\ \  $\forall y_0\in Y, \ \hat{V}(y_0)=\inf_{{\bf u}\in \mathcal{U}} \limsup_{t\to \infty} \gamma_{t}(y_0, {\bf u})$. 
\end{definition}

\begin{proposition}\label{prop:hat-V}
Let $\mathcal{J}=\langle U,f,g\rangle$ be a compact nonexpansive control problems. Then $\hat{V}$ is identical with the general uniform value $V^*$.
\end{proposition}

\textbf{Proof}:   $V^*$ is the general uniform value of $\mathcal{J}$, so $\lim_{t\to\infty} V_t(y_0)=V^*(y_0), \forall y_0\in Y$. On one hand, $V_t(y_0)=\inf_{{\bf u}\in \mathcal{U}}\gamma_{t}(y_0, {\bf u})$ implies that $V^*(y_0)\leq \hat{V}(y_0)$. On the other hand, the existence of the general uniform value $V^*$ implies that: \ $\forall \varepsilon>0, \exists T_0>0, \exists {\bf u},\  s.t. \ \gamma_{t}(y_0, {\bf u})\leq V^*(y_0)+\varepsilon, \forall y_0\in Y, \forall t\geq T_0$.
Thus for all $\varepsilon>0$, $\hat{V}(y_0)\leq \limsup_{t\to\infty} \gamma_t(y_0, {\bf u})\leq V^*(y_0)+\varepsilon$. This completes the proof of the proposition. $Q.E.D.$
 
\begin{definition}\label{def:weighted-average} Let $\beta\in(0,1)$ and $\rho\in(0,1]$, the \textbf{$(\beta,\rho)$-weighted average value} of the control problem $\mathcal{J}=\langle U,f,g\rangle$ is defined as: $\forall y_0\in \RR^d$,
$$\bar{V}_{\rho,\beta}(y_0)=\inf_{\bf u\in \mathcal{U}}\bar{\gamma}_{\rho,\beta}(y_0,{\bf u}), \text{where } \bar{\gamma}_{\rho,\beta}(y_0,{\bf u})=\bigg\{\beta  \gamma_\rho(y_0, {\bf u})+ (1-\beta) \limsup_{t\to\infty} \gamma_t(y_0, {\bf u})\bigg\}.$$
\end{definition}

\begin{proposition}\label{prop:weighted-V} Let $\mathcal{J}=\langle U,f,g\rangle$ be a compact nonexpansive control problem. Then for any $\beta\in(0, 1), \ \lim_{\rho\to 0} \bar{V}_{\rho,\beta}= V^*$.
\end{proposition}

\textbf{Proof}: We fix $\beta\in(0,1)$ and $y_0\in Y$. For each $\rho\in(0,1]$, we use the definition of $\bar{V}_{\rho,\beta}(y_0)$ and Proposition \ref{prop:hat-V} to obtain that: 
 $$\bar{V}_{\rho,\beta}(y_0)=\inf_{\bf u\in \mathcal{U}}\bigg\{\beta  \gamma_\rho(y_0, {\bf u})+ (1-\beta) \limsup_{t\to\infty} \gamma_t(y_0, {\bf u})\bigg\}\geq \beta V_\rho(y_0)+(1-\beta)V^*(y_0).$$  
Thus by the fact that $\lim_{\rho\to 0} V_\rho(y_0)= V^*(y_0)$, we have: $$\liminf_{\rho\to 0}\bar{V}_{\rho,\beta}(y_0)\geq \beta\liminf_{\rho\to 0} V_\rho(y_0)+(1-\beta)V^*(y_0)=V^*(y_0).$$   
Now we prove the other direction $\liminf_{\rho\to 0}\bar{V}_{\rho,\beta}(y_0)\leq V^*(y_0)$. For each pair $\rho\in(0, 1]$ and $t>0$, we consider the $\beta$-weighted evaluation $\theta(\rho, t)$ whose density function is:
$$f_{\theta(\rho, t)}(s)=\beta f_{\theta_\rho}(s)+(1-\beta)f_{\bar{\theta}_t}(s)=\beta \rho e^{-\rho s}+(1-\beta)\frac{1}{t}1_{[0, t]}(s), \forall s\geq 0.$$
$V^*$ is the general uniform value of $\mathcal{J}$, thus for any $\varepsilon>0$, there is some $\eta>0$, $S>0$ and control ${\bf u^*}$ such that: $\gamma_\theta(y_0, {\bf u^*})\leq V^*(y_0)+\varepsilon, \ \forall \theta\in\De(\RR_+)\ \text{ with } \overline{TV}_S(\theta)\leq \eta$. Consider now the evaluation $\theta$ to be separately $(\theta_\rho)$ and $(\bar{\theta}_t)$. An easy computation gives us (cf. Section 3.3 in \cite{Li_2016} for details): 
$$\forall \rho\in(0, 1],  t\geq 0, \ \  \overline{TV}_S(\theta_\rho)\leq 1-e^{-S\rho} \ \text{ and } \    \overline{TV}_S(\bar{\theta}_t)\leq \min\{S/t , 1\}.$$ 
This means that there is $\rho_0>0$ and $T_0>0$ such that: $\forall \rho\leq \rho_0, t\geq T_0$,
$$\overline{TV}_S(\theta_\rho)\leq \eta \text{ and } \overline{TV}_S(\bar{\theta}_t)\leq \eta, \text{ thus }\overline{TV}_S\big(\theta(\rho,t)\big)\leq \eta .$$ 
Taking the convex combination, we obtain that: $\forall (\rho,t)\in(0,\rho_0]\times [T_0,+\infty)$,
\begin{equation}\label{eq:limit_weighted}
\beta \gamma_\rho(y_0, {\bf u^*})+ (1-\beta)\gamma_t(y_0, {\bf u^*})=\gamma_{\theta(\rho,t)}(y_0, {\bf u^*})\leq V^*(y_0)+\varepsilon.
\end{equation}
We take $\limsup$ as $t$ tends to infinity and use the definition of $\bar{V}_{\rho, \beta}$ to obtain:
$$ \bar{V}_{\rho, \beta}(y_0)\leq \beta \gamma_\rho(y_0, {\bf u^*})+ (1-\beta)\limsup_{t\to\infty}\gamma_t(y_0, {\bf u^*})\leq V^*(y_0)+\varepsilon.$$ 
Finally, taking $\liminf$ as $\rho$ tends to zero on both sides of the above inequality yields:
$$\liminf_{\rho\to 0}\bar{V}_{\rho, \beta}(y_0)\leq V^*(y_0)+\varepsilon.$$ 
Since $\varepsilon>0$ can be arbitrarily, the proof is then complete.  $Q.E.D.$

\begin{remarque} 
Equation (\ref{eq:limit_weighted}) implies that the $\varepsilon$-optimal control for the general uniform value $V^*$ is also $\varepsilon$-optimal for all $(\rho, \beta)$-weighted average problems provided that $\rho$ is small enough. In particular, the choice of the weight parameter $\beta$ is irrelevant.     
\end{remarque}

Note that our proof of the above two results replies only on the existence of the general uniform value $V^*$ but not on the particular properties of the compact nonexpansive conditions. 

\section{Proof of Theorem \ref{thm:guv}}\label{sec:proof}
 
\subsection{Some preliminary results}

We introduce several further notations concerning random controls. Our use of random control/strategy in continuous-time control/game problem follows \cite{Cardaliaguet_2007} (see \cite{Aumann_64} for the introduction of randomized strategies for infinite games). We are going to work on a set $\mathcal{S}$ of probability spaces, which is stable by countable product. To fix the ideas, choose $$\mathcal{S}=\{\left([0,1]^n,\mathcal{B}([0,1]^n),\lambda^n\right), \text{ for some } n\in\NN^*\cup\{\infty\}\},$$
where $\mathcal{B}([0,1]^n)$ is the class of Borel sets of $[0,1]^n$, $\lambda^n$ is the Lebesgue measure on $\RR^n$ for $n<\infty$ and $\mathcal{B}([0,1]^\infty)$ is the product Borel-field, $\lambda^\infty$ is the product measure for $n=\infty$.

We might write simply ${\bf{u}}\in\tilde{\mathcal{U}}$ without explicitly mentioning the underlying probability space.  

For any $(y_0,{\bf{u}})\in Y\times \tilde{\mathcal{U}}$ and $t\geq 0$, we denote $\tilde{{\rm \textbf{y}}}(t, {\bf{u}},y_0)=\int_{\Omega}\de_{{\rm \textbf{y}}(t,{\bf u}_{\omega}(y_0), y_0)}d\lambda(\omega)$ for the distribution of the state at time $t$. Any $z\in\De(Y)$ is a probability distribution over $Y$. Let $z$ be the distribution of the initial state,  then $t\mapsto\tilde{{\rm \textbf{y}}}(t, {\bf{u}},z)=\int_{ Y}\int_{\Omega}\de_{{\rm \textbf{y}}(t,{\bf{u}}_{\omega}(y_0), y_0)}d\lambda(\omega)d z(y_0)$ is the expected trajectory in $\De(Y)$ with respect to $z$. The above notations are consistent when a point $y_0\in Y$ is identified with the Dirac measure $\de_{y_0}\in\De(Y)$.

For any $\tilde{\rm\bold{y}}\in\De(Y)$, let $g(\tilde{\rm\bold{y}})=\int_{p\in Y}g(p)d \tilde{\rm\bold{y}}(p)$ be the affine extension of $g$. Using Fubini's theorem, we have: for any $y_0\in Y$ and ${\bf{u}}\in\tilde{\mathcal{U}}$, $$\gamma_\theta(y_0, {\bf{u}})=\int_{[0,+\infty)} \int_{\Omega}g\Big({\rm \textbf{y}}(t, {\bf{u}}_\omega(y_0),y_0)\Big) d\lambda(\omega)d \theta(t)=\int_{[0,+\infty)} g\Big(\tilde{{\rm \textbf{y}}}(t,{\bf{u}},y_0)\Big)d \theta(t).$$

\bigskip 

Before proceeding to the proof of Theorem \ref{thm:guv}, we establish in the rest of this subsection several preliminary results concerning the nonexpansive property of the dynamic in $\De(Y)$.

Let $d_{KR}$ be the \textit{Kantorovich-Rubinstein distance} on $\De(Y)$: $$\forall z, z'\in \De(Y), \ \ d_{KR}(z,z')=\sup_{h\in Lip(1)} \left|\int_Y h d z - \int_Y h d z' \right|,$$
where $Lip(1)$ denotes the set of bounded 1-Lipschitz functions defined on $Y$. 

\begin{lemme}\label{lem:d_KR} For any $z_1, z_2$ in $\De(Y)$ and ${\bf{u}}:Y\times\Omega\times[0,+\infty)\to U$ a random control, there exists some random control ${\bf{v}}:Y\times\hat{\Omega}\times\Omega\times [0,+\infty)\to U$ such that $$d_{KR}\Big(\tilde{{\rm \bold{y}}}(t, {\bf{u}}, z_1), \tilde{{\rm \bold{y}}}(t,{\bf{v}},z_2)\Big)\leq d_{KR}(z_1,z_2), \ \forall t\geq 0.$$
\end{lemme}

 \textbf{Proof}:  We fix ${\bf{u}}:Y\times\Omega\times[0,+\infty)\to U$ a random control defined on the probability space $\big(\Omega, \mathcal{B}(\Omega),\lambda\big)$. We first show the result for $z_1$ and $z_2$ being Dirac measures. Let $p, q$ in $Y$. ${\bf{u}}_\omega(p):[0,+\infty)\to U$ is a pure control for any $\omega\in \Omega$. By the nonexpansive assumption, Proposition 3.6 in \cite{Quincampoix_2011} implies that: for given $\omega\in \Omega$, there is some pure control ${\bf{u}}^{p}(q,\omega,\cdot):[0,+\infty)\to U$, which we denote by ${\bf{u}}^{p}_\omega(q):={\bf{u}}^{p}(q,\omega,\cdot)$, such that
\begin{align}\label{eq:nonexpan-dirac}\Big\|{\rm \bold{y}}\big(t, {\bf{u}}_\omega(p),p\big)-{\rm \bold{y}}\big(t, {\bf{u}}^{p}_\omega(q), q\big)\Big\|\leq \big\|p-q\big\|, \ \forall t\geq 0.\end{align}
Moreover, one can take the mapping $(q,t)\mapsto {\bf{u}}^p_\omega(q,t):={\bf{u}}^{p}(q,\omega,t)$ jointly measurable\footnote{Indeed, consider for (\ref{eq:non-expansive}): we fix any $y_2\in Y$ and $a^\cdot:x\mapsto a^{x}\in U$ measurable. Then we can apply the measurable selection theorem (cf. Theorem 9.1 in \cite{Wagner_77}) for the optimization problem $\inf_{b\in U}\big\langle y_1-y_2, f(y_1,a^{y_1})-f(y_2,b)\big\rangle$ to obtain the existence of a measurable mapping $b^\cdot:x\mapsto b^{x}\in U$ satisfying $\big\langle y_1-y_2, f(y_1,a^{y_1})-f(y_2,b^{y_1})\big\rangle\leq 0, \ \forall y_1\in Y$. The same augument as in the proof of Proposition 3.6 in \cite{Quincampoix_2011} implies the result.}. Letting ${\bf{u}}^p(q,t)=\big({\bf{u}}^p_\omega(q,t)\big)_{\omega\in\Omega}$ for all $t\geq 0$, this defines
a random control:
\begin{align*}
{\bf{u}}^p:Y\times\Omega \times [0,+\infty) & \mapsto U\\
(q, \omega, t)& \mapsto {\bf{u}}^p_\omega(q,t).
\end{align*} From Equation (\ref{eq:nonexpan-dirac}), we deduce that: for any $t\geq 0$,
\begin{equation}\label{eq:u-p_q}
\begin{aligned}
d_{KR}\Big(\tilde{{\rm \textbf{y}}}(t,{\bf{u}},p), \tilde{{\rm \textbf{y}}}(t,{\bf{u}}^p,q)\Big)&= d_{KR}\Big(\int_{\Omega}\de_{{\rm \textbf{y}}(t, {\bf{u}}_\omega(p),p)}d\lambda(\omega), \int_{\Omega}\de_{{\rm \textbf{y}}(t, {\bf{u}}_\omega^p(q),q)}d\lambda(\omega)\Big)\\
&\leq \int_\Omega \Big\|{\rm \textbf{y}}\big(t, {\bf{u}}_\omega(p),p\big)-{\rm \textbf{y}}\big(t, {\bf{u}}_\omega^p(q),q\big)\Big\|d\lambda(\omega)\\
&\leq ||p-q||.
\end{aligned}
\end{equation}
Consider now $z_1, z_2$ in $\De(Y)$. By the Kantorovich-Rubinstein duality formula (cf. Theorem 5.10 in \cite{Villani_2009}), there is a coupling $\xi(\cdot, \cdot)\in\De(Y\times Y)$ with first marginal $z_1$ and second marginal $z_2$, satisfying:
\begin{align}\label{eq:K-duality} d_{KR}(z_1,z_2)=\int_{Y\times Y}\big\|p-q\big\|d\xi(p,q).
 \end{align}
Let $q\mapsto \xi(\cdot|q)\in\De(Y)$ be the conditional distribution of $\xi$ on its first marginal. Consider now the Borel isomorphic mapping (between two standard Borel spaces)
$$h_q: \hat{\Omega}=\Big([0,1], \mathcal{B}\big([0,1]\big), \lambda\Big) \mapsto \Big( Y,\mathcal{B}(Y), \xi(\cdot|q)\Big)$$
where $$ \forall B\in\mathcal{B}(Y), \ \xi(B|q)=\lambda\left(h_q^{-1}(B)\right).$$ 
Define now
\begin{equation*}
\begin{aligned}
{\bf{v}}: Y\times \hat{\Omega} &\times \Omega\times [0,+\infty) \ \  \ \mapsto \  \ \ U \\
& (q, \hat{\omega}, \omega, t)  \ \ \ \ \  \ \ \ \ \ \mapsto \ \  \ {\bf{v}}_{(\hat{\omega}, \omega)}(q,t)={\bf{u}}^{h_q(\hat{\omega})}_\omega(q, t),
\end{aligned}
\end{equation*}
which is jointly measurable as the composition function of ${\bf{u}}_\omega^\cdot(q,t)$ and $h_q(\hat{\omega})$. $(\Omega_{\bf{v}},{\bf{v}})$ is then a random control defined on the product Borel probability space $\big(\hat{\Omega}\otimes\Omega, \mathcal{B}(\hat{\Omega}\otimes\Omega), \lambda^2\big)$. The interpretation of the random control ${\bf{v}}$ is that: at each initial point $q\in Y$, ${\bf{v}}(q, \cdot)$ randomly takes the value ${\bf{u}}^p(q, \cdot)$ according to the probability law $d\xi( p|q)$.

Next, we check that ${\bf{v}}$ satisfies the nonexpansive condition. For any $t\geq 0$ and $q\in supp (z_2)$, we use Fubini's theorem (first on $\hat{\Omega}\times \Omega$ and then on $Y\times Y$) and the change of variable "$p=h_q(\hat{\omega})$" to obtain:
\begin{align*}
\tilde{{\rm \textbf{y}}}(t, {\bf{v}}, z_2)=\int_{Y} \tilde{{\rm \textbf{y}}}\big(t, {\bf{v}},q\big)d z_2(q)=&\int_{Y}  \int_{\hat{\Omega}\times \Omega}\de_{{\rm \textbf{y}}\big(t, {\bf{u}}^{h_q(\hat{\omega})}_\omega(q),q\big)} d \lambda^2(\hat{\omega},\omega)d z_2(q)\\
=&\int_{Y}  \int_{\hat{\Omega}}\Bigg[\int_{ \Omega}\de_{{\rm \textbf{y}}\big(t, {\bf{u}}^{h_q(\hat{\omega})}_\omega(q),q\big)} d \lambda(\hat{\omega})\Bigg]d \lambda(\omega)d z_2(q)\\
=&\int_{Y} \int_Y \Bigg[\int_{\Omega}\de_{{\rm \textbf{y}}\left(t, {\bf{u}}^p_\omega(q),q\right)} d \lambda(\omega)\Bigg])d \xi(p|q)d z_2(q)\\
=&\int_{Y\times Y} \tilde{{\rm \textbf{y}}}\big(t,{\bf{u}}^p,q\big) d \xi(p,q).
\end{align*}
We then deduce that:
\begin{align*}
d_{KR}\Big(\tilde{{\rm \textbf{y}}}(t, {\bf{u}}, z_1),\tilde{{\rm \textbf{y}}}(t, {\bf{v}}, z_2) \Big)&=d_{KR}\Big(\int_{Y}\tilde{{\rm \textbf{y}}}\big(t, {\bf{u}}, p\big)d z_1(p), \ \int_{Y\times Y} \tilde{{\rm \textbf{y}}}\big(t,{\bf{u}}^p,q\big) d \xi(p,q)\Big)\\
&= d_{KR}\big(\int_{Y\times Y}\tilde{{\rm \textbf{y}}}(t, {\bf{u}}, p)d\xi(p,q),\  \int_{Y\times Y} \tilde{{\rm \textbf{y}}}\big(t,{\bf{u}}^p,q\big) d \xi(p,q)\big)\\
&\leq \int_{Y\times Y} d_{KR}\Big(\tilde{{\rm \textbf{y}}}\big(t, {\bf{u}}, p\big), \tilde{{\rm \textbf{y}}}\big(t,{\bf{u}}^p,q\big) \Big)d\xi(p,q)\\
\Big(\text{by } (\ref{eq:u-p_q}) \Big)& \leq \int_{Y\times Y} ||p-q||d \xi(p,q)\\
\Big(\text{by } (\ref{eq:K-duality}) \Big) & = d_{KR}(z_1,z_2).
\end{align*}
This completes our proof for the lemma.   \  $Q.E.D.$  
 
\bigskip 

The proof of Theorem \ref{thm:guv} involves some compact property of the set of random controls for compact non-expansive control problem, as we will show in Lemma \ref{lem:limit-traj}:  the "limit trajectory" in $\De(Y)$ (cf. Lemma \ref{def:limit-trajectory}) of a sequence of random controls can be arbitrarily approximated by the trajectory induced by one random control. Here the random control that we are going to construct is defined as concatenations of certain sequence of random controls, namely \textit{behavior control}. Formal definitions are as follows. 

\begin{definition} For any two random controls $(\Omega_{\bf{u}},{\bf{u}})$ and $(\Omega_{\bf{v}},{\bf{v}})$ and a time $T>0$. The \textbf{concatenation} of ${\bf{u}}$ and ${\bf{v}}$ at time $T$ is defined as the random control $(\Omega_{\bf{u}}\times \Omega_{\bf{v}}, {\bf{u}}\oplus_{T}{\bf{v}})$ with: $\forall \big(y_0, (\omega_1,\omega_2),t\big)\in Y\times (\Omega_{\bf{u}}\times \Omega_{\bf{v}})\times [0,+\infty)$,
$$[{\bf{u}}\oplus_{T}{\bf{v}}] _{(\omega_1,\omega_2)}(y_0,t)=1_{\{t<T\}}{\bf{u}}_{\omega_1}(y_0,t)+1_{\{t\geq T\}}{\bf{v}}_{\omega_2}\Big({\rm \bold{y}}\big(y_0,{\bf{u}}_{\omega_1}(y_0),T\big),t-T\Big).$$
\end{definition} 

\begin{definition}\label{def:behavior control} Fix $0=t_0<t_1\cdot\cdot\cdot< t_m<\cdot\cdot\cdot$ a partition of $[0,+\infty)$, and $\left(\Omega_m,{\bf{u}}^m\right)_{m\geq 1}$ a sequence of random controls. Let $\big(\otimes_{m'=1}^m \Omega_{m'}, {\bf{u}}^{[m]}\big)_{m\geq 1}$ be a sequence of random controls defined inductively as: ${\bf{u}}^{[1]}={\bf{u}}^1$ and ${\bf{u}}^{[m+1]}={\bf{u}}^{[m]}\oplus_{t_{m}} {\bf{u}}^{m+1}$ for any $m\geq 1$. The \textbf{behavior control} $$\big(\Omega_{[\infty]},{\bf{u}}^{[\infty]}\big):=\big(\otimes_{m\geq 1}\Omega_m, {\bf{u}}^1\oplus_{t_1}\cdot\cdot\cdot {\bf{u}}^t\oplus_{t_m}\cdot\cdot\cdot\big)$$ is defined as the concatenations of $({\bf{u}}^m)$ sequentially at points $(t_m)$:\\
 for any $\big(y_0, (\omega_m)_{m\geq 1},t\big)$ in $Y\times \Omega_{[\infty]}\times [0,+\infty)$, 
$${\bf{u}}^{[\infty]}_{(\omega_m)_{m\geq 1}}(y_0,t)=\sum_{m\geq 1}1_{\{t_{m-1}\leq t<t_m\}}{\bf{u}}^{[m]}_{\omega^m}(y_0,t), \text{\ where \ }   \omega^m:=(\omega_1,...,\omega_m).$$
\end{definition}
 
\begin{remarque} \label{rem: behavior} The behavior control $(\Omega_{[\infty]}, {\bf{u}}^{[\infty]})$ is also a random control with the product Borel space $\Omega_{[\infty]}= \otimes_{m\geq 1}\Omega_m$, which, as a countable union, is still in $\mathcal{S}$. 
\end{remarque}

\begin{remarque} \label{rem: kuhn}
On the other hand, from Kuhn's theorem (cf. \cite{Aumann_64}, Section 5): for any random control, one can construct a behavior control such that the trajectories in $\De(Y)$ induced by them are the same. More precisely, we fix $(t_m)_{m\geq 0}$ a  partition of $[0,+\infty)$, $\big(\Omega_{\bf{u}}, {\bf{u}}\big)$ a random control and $y_0\in Y$ an initial state. Then, there exists some behavior control $\big(\otimes_{m\geq 1}\Omega_m,\overline{{\bf{u}}}\big)$ as concatenations of some sequence of random controls $\big(\Omega_m,\overline{{\bf{u}}}^m\big)_{m\geq 1}$ at points $(t_m)$ such that starting from $y_0$, the trajectories in $\De(Y)$ generated by both ${\bf u}$ and $\overline{{\bf{u}}}$ are the same, i.e. $\tilde{{\rm \bold{y}}}(t, {\bf{u}},y_0)=\tilde{{\rm \bold{y}}}(t, \overline{{\bf{u}}}, y_0)$,  $\forall t\geq 0, a.e.$
\end{remarque}
  
Fix any $y_0$ in $Y$. $\big(\tilde{{\rm \bold{y}}}(\cdot,{\bf{u}}^k, y_0)\big)_{k\geq 1}$ is the sequence of trajectories in $\De(Y)$ generated by a sequence of random controls $({\bf{u}}^k)_{k\geq 1}$ with the same initial point $y_0$.

\begin{definition}\label{def:limit-trajectory}  A measurable mapping $t\mapsto \overline{\bf{y}}(t)$ from $\RR_+$ to $\big(\De(Y),d_{KR}\big)$ is a \textbf{limit trajectory} of the sequence $\big(\tilde{{\rm \bold{y}}}(\cdot,{\bf{u}}^k, y_0)\big)_{k\geq 1}$ if there is some $\psi(k)$ such that for any $m\geq 1$, $\tilde{{\rm \bold{y}}}(\cdot,{\bf{u}}^{\psi(k)}, y_0)$ converges (for $d_{KR}$) to $\overline{\bf{y}}(\cdot)$ uniformly on $[0,m]$ as $k$ tends to infinity. 
\end{definition}
 
We first show that the limit trajectory exists for any sequence. 

\begin{lemme}\label{lem:limit-traj}  \  $\big(\tilde{{\rm \bold{y}}}(\cdot,{\bf{u}}^k, y_0)\big)_{k\geq 1}$ has a limit trajectory in $\De(Y)$.  
\end{lemme}

\textbf{Proof}:  \  Fix an $m\geq 0$, we look at the restriction of each $\tilde{{\rm \bold{y}}}(\cdot,{\bf{u}}^k, y_0)$ on the compact interval $[m, m+1]$. Then the family $\{\tilde{{\rm \bold{y}}}(\cdot,{\bf{u}}^k, y_0):k\geq 1\}$ are continuous mappings from $[m,m+1]$ to the compact domain $\big(\De(Y),d_{KR}\big)$. One can use Ascoli's theorem to deduce the existence of a uniform convergent subsequence $\big(\tilde{{\rm \bold{y}}}(\cdot,{\bf{u}}^{\psi(k)}, y_0)\big)_{k\geq 1}$ on $[m, m+1]$. To obtain this, it is sufficient for us to prove that the family $\big\{\tilde{{\rm \bold{y}}}(\cdot,{\bf{u}}^k, y_0):k\geq 1\big\}$ \big(restricted on $[m, m+1]$\big) is equicontinuous. 

We fix any $k\geq 1$ and $s,t\in[m, m+1]$. Then by the definition of $d_{KR}$, we deduce that
\begin{equation}\label{eq:ascoli}
\begin{aligned}
d_{KR}\Big(\tilde{{\rm \bold{y}}}(t,{\bf{u}}^k, y_0), \tilde{{\rm \bold{y}}}(s,{\bf{u}}^{k}, y_0)\Big)& \leq \int_{\Omega}\Big\|{\rm \bold{y}}\big(t,{\bf{u}}^k_\omega(y_0), y_0\big)-{\rm \bold{y}}\big(s,{\bf{u}}^k_\omega(y_0), y_0\big)\Big\| d\lambda(\omega)\\
 & \leq a\big(1+\sup_{y\in Y}\|y\|\big)  |t-s| ,
\end{aligned}
\end{equation}
where we have used in the last inequality the fact that $t\mapsto {\rm\bold{y}}\big(t,{\bf{u}}^k_\omega(y_0), y_0\big)$ is absolutely continuous \big(cf. assumptions in ($\ref{ex:assumptions}$) \big) . As $Y\subseteq \RR^d$ is compact, Equation (\ref{eq:ascoli}) proves that the family  $\big\{\tilde{{\rm \bold{y}}}(\cdot,{\bf{u}}^k, y_0):k\geq 1\big\}$ \big(restricted on $[m, m+1]$\big) is equicontinuous. By extracting subsequences for each $m$, we obtain the existence of a limit trajectory of $\big(\tilde{{\rm \bold{y}}}(\cdot,{\bf{u}}^k, y_0)\big)_{k\geq 1}$.  \  $Q.E.D.$

\begin{lemme} \label{lem:limit-nonexpan} \  Let $\overline{{\rm \bold{y}}}(\cdot):t\mapsto \overline{{\rm \bold{y}}}(t)$ be a limit trajectory of $\big(\tilde{{\rm \bold{y}}}(\cdot,{\bf{u}}^k,y_0)\big)_k$. Then for any $\varepsilon>0$, there is some behavior control ${\bf{u}}^*$ whose trajectory in $\De(Y)$ is $\varepsilon$-close to $\overline{\bf{y}}(\cdot)$ along time, i.e.,
$$\forall \varepsilon>0, \ \exists {\bf{u}}^*\in\tilde{\mathcal{U}},\ s.t.\ d_{KR}\Big(\tilde{{\rm \bold{y}}}(t,{\bf{u}}^*,y_0), \overline{\rm \bold{y}}(t)\Big)\leq \varepsilon, \ \  \forall t\geq 0.$$
\end{lemme}

\textbf{Proof}: \ The idea is to construct a behavior control ${\bf{u}}^*$ by consecutive  intervals, such that on each of them, ${\bf{u}}^*$ follows one random control in the family $\{{\bf u}^k\}$  whose trajectory is close to the limit $\overline{\bf{y}}$. The proof relies on the nonexpansive property established in Lemma \ref{lem:d_KR}, which ensures that by iteration, the trajectory generated by ${\bf{u}}^*$ is close to $\overline{\bf{y}}$ on the whole $\RR_+$. 

Let $\varepsilon>0$ be fixed. The behavior control ${\bf{u}}^*$ will be constructed as concatenations of a sequence of random controls $(\hat{\bf{u}}^{K_m})$ (to be specified later on) at points $\{1,2,3,...\}$.
 
By definition, $t\mapsto \overline{\rm \bold{y}}(t)$ is a limit trajectory of $\big(\tilde{{\rm \bold{y}}}(\cdot, {\bf{u}}^k, y_0)\big)_k$ in $\De(Y)$ for $d_{KR}$, so for each $m\geq 0$, there exists some $K_{m+1}>0$ such that:
\begin{align}\label{eq:non-expan-1}
d_{KR}\Big(\tilde{{\rm \bold{y}}}(t, {\bf{u}}^{K_{m+1}},  y_0), \overline{\rm \bold{y}}(t)\Big)\leq \varepsilon^{m+1}, \ \forall t\in[m, m+1].
\end{align}

Following Remark \ref{rem: kuhn}, we could have assumed that each ${\bf{u}}^{K_{m+1}}$ is a behavior control, and let $\overline{\bf{u}}^{K_{m+1}}: Y\times \Omega\times [0,+\infty)\to U$ be the component of ${\bf{u}}^{K_{m+1}}$ on interval $[m, m+1]$.

In order to define the behavior control ${\bf{u}}^*$ , it is sufficient to construct by induction a sequence of random controls $(\hat{\bf{u}}^{K_{m}})_{m\geq 1}$ such that for all $m\geq 1$, $$d_{KR}\Big(\tilde{{\rm \bold{y}}}(t,{\bf{u}}^{[m]},y_0), \overline{\rm \bold{y}}(t)\Big)\leq 2\sum_{\ell=1}^m\varepsilon^\ell, \ \  \forall t\in [0,m],$$ where $${\bf{u}}^{[1]}= \hat{\bf{u}}^{K_1} \text{ \ and \ } {\bf{u}}^{[m]}= \hat{\bf{u}}^{K_{1}}\oplus_{1}\cdot\cdot\cdot\oplus_{m-1}\hat{\bf{u}}^{K_{m}}, \ \ m\geq 2.$$

For $m=1$, let $\hat{\bf{u}}^{K_1}=\overline{\bf{u}}^{K_1}$, then from Equation (\ref{eq:non-expan-1}), $d_{KR}\Big(\tilde{{\rm \bold{y}}}(t_1, {\bf{u}}^{[1]}, y_0),{\bf \overline{y} } (t)\Big)\leq \varepsilon, \ \forall t\in[0,1]$. This initializes our induction.

Assume that $\hat{\bf{u}}^{K_1},...,\hat{\bf{u}}^{K_m}$ are defined and let ${\bf{u}}^{[m]}=\hat{\bf{u}}^{K_{1}}\oplus_{1}\cdot\cdot\cdot\oplus_{m-1}\hat{\bf{u}}^{K_{m}}$ (${\bf{u}}^{[1]}=\hat{\bf{u}}^{K_1}$ for $m=1$) satisfy:
\begin{align}\label{eq:non-expan-2}
d_{KR}\Big(\tilde{{\rm \bold{y}}}(t,{\bf{u}}^{[m]},y_0), \overline{\rm \bold{y}}(t)\Big)\leq 2\sum_{l=1}^m\varepsilon^{l}, \ \forall t\in[0,m].
\end{align}
Next we construct the control $\hat{\bf{u}}^{K_{m+1}}$ thus complete the definition of ${\bf{u}}^{[m+1]}$ on $[m, m+1)$. To do this, we consider the two distributions $\tilde{{\rm \bold{y}}}(m,{\bf{u}}^{[m]}, y_0)$ and $\tilde{{\rm \bold{y}}}(m, {\bf{u}}^{K_{m+1}}, y_0)$. Take $t=m$ in Equation (\ref{eq:non-expan-1}) and in Equation (\ref{eq:non-expan-2}), we use the triangle inequality obtain a bound on the distance between them:
\begin{equation}\label{eq:non-expan-3}
\begin{aligned}
& \ d_{KR}\Big(\tilde{{\rm \bold{y}}}(m,{\bf{u}}^{[m]}, y_0), \tilde{{\rm \bold{y}}}(m, {\bf{u}}^{K_{m+1}}, y_0)\Big) \leq & \ 2\sum_{l=1}^m\varepsilon^{l}+\varepsilon^{m+1}.
\end{aligned}
\end{equation}
Consider then the random control $\overline{\bf{u}}^{K_{m+1}}$ on the initial distribution $\tilde{{\rm \bold{y}}}(m, {\bf{u}}^{K_{m+1}}, y_0)$, and apply Lemma \ref{lem:d_KR} to deduce the existence of some random control $\hat{\bf{u}}^{K_{m+1}}$ on the starting distribution $\tilde{\rm\bold{y}}(m,{\bf{u}}^{[m]},y_0)$ such that:
\begin{equation}\label{eq:nin-expan-4}
\begin{aligned}
& \ d_{KR}\Big(\tilde{{\rm \bold{y}}}\big(\De,\hat{\bf{u}}^{K_{m+1}}, \tilde{{\rm \bold{y}}}(m,{\bf{u}}^{[m]}, y_0)\big), \ \tilde{{\rm \bold{y}}}(\De,\bar{\bf{u}}^{K_{m+1}},\tilde{{\rm \bold{y}}}(m, {\bf{u}}^{K_{m+1}}, y_0)\Big)\\
\leq & \ d_{KR} \Big( \tilde{{\rm \bold{y}}}(m,{\bf{u}}^{[m]}, y_0),  \  \tilde{{\rm \bold{y}}}(m, {\bf{u}}^{K_{m+1}}, y_0)\Big)\\
\leq & \ 2 \sum_{l=1}^{m}\varepsilon^l+\varepsilon^{m+1}, \ \forall \De\in[0, 1].
\end{aligned}
\end{equation}

By definition of $\overline{\bf{u}}^{K_{m+1}}$, we have that for all $\De\in[0,1]$,
\begin{align}\label{eq:non-expan-5}
\tilde{{\rm \bold{y}}}\big(\De,\overline{\bf{u}}^{K_{m+1}},\tilde{{\rm \bold{y}}}(m, {\bf{u}}^{K_{m+1}}, y_0)\big)=\tilde{{\rm \bold{y}}}\big(\De+m, {\bf{u}}^{K_{m+1}}, y_0\big).
\end{align}
Define now $\hat{\bf{u}}^{[m+1]}=\hat{\bf{u}}^{[m]}\oplus_{m}\hat{\bf{u}}^{K_m+1}$. This gives us:
\begin{align}\label{eq:non-expan-6}
\tilde{{\rm \bold{y}}}\big(\De,\hat{\bf{u}}^{K_{m+1}}, \tilde{{\rm \bold{y}}}(m,{\bf{u}}^{[m]}, y_0)\big)=\tilde{{\rm \bold{y}}}\big(\De+m, {\bf{u}}^{[m+1]}, y_0\big).
\end{align}

We substitute Equation (\ref{eq:non-expan-5}) and Equation (\ref{eq:non-expan-6}) back into Equation (\ref{eq:nin-expan-4}) and use the change the variable to obtain that:
\begin{align*}
\forall t\in[m,m+1], \ \  d_{KR}\Big(\tilde{{\rm \bold{y}}}(t, {\bf{u}}^{[m+1]},y_0),\  \tilde{{\rm \bold{y}}}(t, {\bf{u}}^{K_{m+1}},y_0) \Big)\leq 2\sum_{l=1}^{m}\varepsilon^l+\varepsilon^{m+1}.
\end{align*}
Finally, with the help of definition of $K_{m+1}$ in Equation (\ref{eq:non-expan-1}), we obtain that: for any $t\in[m,{m+1}]$,
\begin{align*}
d_{KR}\Big(\tilde{{\rm \bold{y}}}(t, {\bf{u}}^{[m+1]},y_0),\  \bar{{\rm \bold{y}}}(t) \Big)\leq 2\sum_{l=1}^{m}\varepsilon^l+\varepsilon^{m+1}+\varepsilon^{m+1}=2\sum_{l=1}^{m+1}\varepsilon^l.
\end{align*}
This finishes the inductive definition of the sequence $\big(\hat{\bf{u}}^{K_m}\big)_{m\geq 1}$.

To conclude, we set the behavior control  $${\bf{u}}^*=\hat{\bf{u}}^{K_1}\oplus_{1}\cdot\cdot\cdot \hat{\bf{u}}^{K_m}\oplus_{m}\cdot\cdot\cdot,$$
as concatenations of $\big(\hat{\bf{u}}^{K_m}\big)_{m\geq 1}$ at points $\{m\geq 1\}$, and by our inductive construction:
$$\forall t\geq 0, \ d_{KR}\Big(\tilde{\rm\bold{y}}(t,{\bf{u}}^*,y_0) ,\ \bar{\rm \bold{y}}(t)\Big)\leq 2\sum_{\ell=1}^\infty \varepsilon^\ell=\frac{2\varepsilon}{1-\varepsilon}\leq 3\varepsilon,$$
as long as $\varepsilon\in(0,\frac{1}{3}]$. This completes our proof for the lemma by considering $\varepsilon'=\varepsilon/3$. \ $Q.E.D.$

\subsection{Main part of the proof}
\noindent This section is devoted for the proof of Theorem \ref{thm:guv}, which is divided into three parts. 

\textbf{Part A} aims at establishing certain optimality properties for a sequence of controls (Lemma \ref{lem:uniform-1}); In \textbf{Part B}, we use the compact nonexpansive property (Lemma \ref{lem:limit-nonexpan}) to obtain a "limit" control of the above sequence ensuring that the average cost on each (consecutive) block of fixed length is no more than $V^*$ \big(Equation (\ref{eq:limit-close-control})\big); \textbf{Part C} concludes the proof through a comparison of the (normalized) $\theta$-evaluated payoff to the average cost by blocks.   

\bigskip

\noindent\textbf{Part A}.  \ For any $t\geq 0$ , $S>0$, $y_0\in Y$ and ${\bf{u}}\in \tilde{\mathcal{U}}$, denote $$\gamma_{t,S}(y_0,{\bf{u}})=\frac{1}{S}\int_{[t,t+S]}g\big(\tilde{\rm\bold{y}}(s,{\bf{u}},y_0)\big)d s, \ \forall y_0\in Y.$$
\noindent For $T\geq 0$, we put 
\begin{align*}
\varphi_{T,S}(y_0) = \inf_{{\bf{u}}\in\mathcal{\tilde{U}}} \sup_{\mu\in \De([0,T])} \int_{[0,T]}\gamma_{t, S}(y_0,{\bf{u}})d\mu(t).
\end{align*}

Fix any $T,S,y_0$. We first prove a minmax result for $\varphi_{T,S}(y_0)$. We denote for each $s\geq 0$: $$\beta_s(\mu, S)=\frac{1}{S}\int_{\max\{0,\ s-S\}}^{\min\{T,s\}}d\mu(t)=\frac{1}{S}\mu\Big([0,T]\cap [s-S,s]\Big).$$ 
Then from the definition of $\gamma_{t,S}(y_0,{\bf{u}})$, we obtain that
\begin{eqnarray*}
\int_{[0,T]}  \gamma_{t, S}(y_0,{\bf{u}})d\mu(t)&=&\int_{[0,T]} \Bigg(\frac{1}{S}\int_{[t,t+S]} g\big(\tilde{\rm\bold{y}}(s,{\bf{u}},y_0)\big)ds\Bigg)d\mu(t)\\
(''\text{Fubini's theorem}'') &=& \int_{[0,T+S]} \beta_s(\mu, S) g\left(\tilde{\rm\bold{y}}(s,{\bf{u}}, y_0)\right)d s.
\end{eqnarray*}

Note that for each fixed $S>0$ and $\mu\in\De([0,T])$, the mapping $t\mapsto \beta_t(\mu,S)$ defines a density function of some evaluation over $\RR_+$, which we denote by $\varsigma(\mu,S)$. This enables us to write
$$\varphi_{T,S}(y_0)=\inf_{\bf{u}\in \tilde{\mathcal{U}}} \sup_{\mu\in \De([0,T])} \gamma_{\varsigma(\mu, S)} (y_0,{\bf{u}}).$$

Next, we use Sion's minmax theorem (cf.  Appendix A.3 in \cite{Sorin_2002}) to show that the operators "$\inf$" and "$\sup$" of the above equation commute. Indeed, $\tilde{\mathcal{U}}$ is convex; $\De([0,T])$ is convex and weak-* compact and the payoff function $(\mu, {\bf{u}})\mapsto\gamma_{\varsigma(\mu,S)} (y_0,{\bf{u}})$ is affine in both $\mu$ and ${\bf{u}}$; moreover the function $\gamma_{t,S}(y_0,{\bf{u}})$ is continuous in $t$ for given ${\bf{u}}$ ($g$ is continuous in $y$ and each trajectory is absolutely continuous), and so is ${\bf{u}}\mapsto\gamma_{\varsigma(\mu,S) } (y_0,{\bf u} )$. Then we obtain:
\begin{equation}\label{eq:minmax}
\begin{aligned}\varphi_{T,S}(y_0)
&=\sup_{\mu\in \De([0,T])} \inf_{{\bf{u}} \in \tilde{\mathcal{U}}} \gamma_{\varsigma(\mu, S)}(y_0,{\bf{u}})=\sup_{\mu\in \De([0,T])} V_{\varsigma(\mu, S)}(y_0).
\end{aligned}
\end{equation} 
 
\begin{lemme}\label{lem:uniform-1}
For any $\varepsilon>0$, there is some $S_0> 0$ such that
$$\forall T\geq 0, \exists {\bf{u}}^T\in\tilde{\mathcal{U}}: \forall y_0\in Y, \ \gamma_{t,S_0}(y_0, {\bf{u}}^T)\leq V^*(y_0)+\varepsilon,\forall t\leq T.$$
\end{lemme}

\textbf{Proof}:  Let $\theta\in\De(\RR_+)$ be an evaluation that is absolutely continuous with respect to the Lebesgue measure, thus it admits a density function $t\mapsto f_\theta(t)$. For any $s\geq 0$, we denote $I_s(\theta)=\int_{[0,+\infty)} \big|f_\theta(t+s) - f_\theta(t)\big|d t$.  Then one obtains that (cf. Lemma 3.3 in \cite{Li_2016}):
$$I_s(\theta)/2 \leq TV_s(\theta)\leq I_s(\theta).$$
For any $S>0$, $T\geq 0$ and $\mu\in\big([0,T]\big)$, we apply the above expression for $\varsigma(\mu,S)$ so as to obtain the following bound: $\forall s\in[0,S]$,
\begin{equation}\label{eq:bound}
\begin{aligned}
& I_s\big(\varsigma(\mu, S)\big) \\
= &\frac{1}{S}\int_{[0,T+S]} \Big[\mu\big([t-S,t-S+s]\cap [0,T]\big)+\mu\big([t,t+s]\cap [0,T]\big)\Big]d t\\
=& \int_{[0,T+S]}\int_{[t-S,t-S+s]\cap [0,T]}d \mu(s')d t+\int_{[0,T+S]}\int_{[t,t+s]\cap [0,T]}d \mu(s')d t\\
=& \frac{s}{S}\cdot \Big(\mu\big([-S,T+s]\cap [0,T]\big) + \mu\big([0,s+T+S]\cap [0,T]\big)\Big) \\
\leq & \frac{2s}{S},
\end{aligned}
\end{equation}
where we have used Fubini's theorem to obtain the third equality.

According to Proposition \ref{prop:LQR}, the general limit value exists and is equal to $V^*$, i.e. for any $\varepsilon>0$, we take $\eta>0$ and $S'>0$ such that:
\begin{align}\label{eq:limit-value}\forall \theta\in\De(\RR_+), \ \  \left(\sup_{0\leq s\leq S'}TV(\theta)\leq \eta' \Longrightarrow \left(\forall y_0\in Y, \ |V_\theta(y_0)-V^*(y_0)|\leq \varepsilon/2\right)\right).
\end{align}

Take $S_0=\max\{\frac{2S'}{\eta'},S'\}$. From Equation (\ref{eq:bound}), we obtain that: 
$$\forall s\in[0, S'], S\geq S_0, T>0, \mu\in\De\big([0,T]\big), \ TV_s\left(\varsigma(\mu, S)\right)\leq \frac{2s}{S}\leq \frac{2S'}{S_0}\leq \eta',$$ thus by Equation (\ref{eq:limit-value}): \ $\forall y_0\in Y, \ |V_{\varsigma(\mu,S)}(y_0)-V^*(y_0)|\leq \varepsilon/2$. 

Next, from Equation (\ref{eq:minmax}),  $\varphi_{T,S}(y_0)=\sup_{\mu\in\De([0,T])} V_{\varsigma(\mu,S)}(y_0)$, we deduce that
\begin{align}\label{eq:beta-theta_n}
\forall \varepsilon>0, \exists S_0> 0:   \forall S\geq S_0, \forall T\geq 0, \forall y_0\in Y, \ |\varphi_{T,S}(y_0)-V^*(y_0)|\leq \varepsilon/2.
\end{align}

Finally, according to the definition $\varphi_{T,S}(y_0)=\inf_{\bf{u}}\sup_{\mu} \gamma_{\varsigma(\mu,S)}(y_0,{\bf{u}})$, we have $$\forall T>0, \exists {\bf{u}}^T\in\tilde{\mathcal{U}}: \forall y_0\in Y,\  \gamma_{t,S_0}(y_0,{\bf{u}}^T)\leq \varphi_{T,S_0}(y_0)+\varepsilon/2, \forall t\in[0,T].$$ 

Together with Equation \ref{eq:beta-theta_n}, one obtains
\begin{align}\label{eq:S_0-t}
\forall \varepsilon>0, \exists S_0> 0: \forall T\geq 0, \exists {\bf{u}}^T\in\tilde{\mathcal{U}}: \forall y_0\in Y, \ \gamma_{t,S_0}(y_0, {\bf{u}}^T)\leq V^*(y_0)+\varepsilon,\forall t\leq T.
\end{align}

The proof of the lemma is then complete. $Q.E.D.$

\bigskip 

\noindent \textbf{Part B}. \  Fix now any $\varepsilon>0$, and consider $S_0>0$ and the random control ${\bf{u}}^T\in \tilde{\mathcal{U}}$ for any $T>0$ given as in Lemma \ref{lem:uniform-1}. We take an increasing sequence $(T_k)_{k\geq 1}$ in $\RR_+$ and fix any $y_0\in Y$. For each $k\geq 1$, $t\mapsto \tilde{\rm\bold{y}}(t,{\bf{u}}^{T_k},y_0)$ is the trajectory of ${\bf{u}}^{T_k}$ in $\De(Y)$. Thus from Equation (\ref{eq:S_0-t}), we obtain: 
\begin{align}\label{eq:sequence} \gamma_{t,S_0}(y_0,{\bf{u}}^{T_k})=\frac{1}{S_0}\int_{[t,t+S_0]}g\left(\tilde{\rm\bold{y}}(t,{\bf{u}}^{T_k},y_0)\right)d s\leq V^*(y_0)+\varepsilon \text{ for all } t\leq T_k. 
\end{align}                                                                                                                                                                                                                                                                                                                                                                                                
                                                                                                                                                                                                                                                                                                                                                                                        Let $\bar{\rm \bold{y}}(\cdot):t\mapsto \bar{\rm \bold{y}}(t)$ be a limit trajectory of the sequence $\big(\tilde{\rm\bold{y}}(\cdot,{\bf{u}}^{T_k},y_0)\big)_{k\geq 1}$  i.e. there is a subsequence $\psi(k)$ such that $\tilde{\rm\bold{y}}(\cdot, {\bf{u}}^{T_{\psi(k)}},y_0)$ converges uniformly to $\bar{\bf{y}}(\cdot)$ on each $[m,{m+1}]$. Since $g$ is continuous on the compact invariant set $Y$, and $\De(Y)$ is weak-* compact for the topology induced by the distance $d_{KR}$ \big(cf. Theorem 6.9 in \cite{Villani_2009}\big), we let $k$ tend to infinity \big(along the subsequence $\psi(k)$\big) in Equation (\ref{eq:sequence}) to get
\begin{equation}\label{eq:limit-close-trajectory}
\frac{1}{S_0}\int_{[t,t+S_0]}g\big(\bar{\rm \bold{y}}(s)\big)d s\leq V^*(y_0)+\varepsilon, \ \forall t\geq 0.
\end{equation}

Now we apply Lemma \ref{lem:limit-nonexpan} for the sequence $\big(\tilde{\rm\bold{y}}(\cdot,{\bf{u}}^{T_k},y_0)\big)_{k\geq 1}$ and its limit trajectory $\bar{\rm \bold{y}}(\cdot)$ to obtain the existence of some behavior control ${\bf{u}}^*$ such that: $$d_{KR}\big(\tilde{\rm\bold{y}}(t,{\bf{u}}^*,y_0),\bar{\rm \bold{y}}(t)\big)\leq \varepsilon, \ \forall t\geq 0.$$
Together with Equation (\ref{eq:limit-close-trajectory}), we obtain:
\begin{equation}\label{eq:limit-close-control}
 \gamma_{t,S_0}(y_0, {\bf{u}}^*)=\frac{1}{S_0}\int_{[t,t+S_0]}g\big(\tilde{\rm\bold{y}}(s,{\bf{u}}^*,y_0)\big)d s\leq V^*(y_0)+2\varepsilon, \ \forall t\geq 0.
\end{equation}

\bigskip

\noindent \textbf{Part C}. \ Let $\theta\in\De(\RR_+)$ be any evaluation. We integrate Equation (\ref{eq:limit-close-control}) over $t\geq 0$ with respect to $\theta$, to obtain:
\begin{equation}\label{eq:uniform-2}
\begin{aligned}
 V^*(y_0)+2\varepsilon \geq & \int_{[0,+\infty)}  \gamma_{t, S_0}(y_0,{\bf{u}}^*)d\theta(t)\\ 
 = &\int_{[0,+\infty)} \left(\frac{1}{S_0}\int_{[t,t+S_0]} g\left(\tilde{\rm\bold{y}}(s,{\bf{u}}^*,y_0)\right)d s\right)d\theta(t)\\
(''\text{Fubini's theorem}'') =& \int_{[0,+\infty)} \beta_s(\theta, S_0) g\left(\tilde{\rm\bold{y}}(s,{\bf{u}}^*,y_0)\right)d s\\
=& \gamma_{\varsigma(\theta, S_0)}(y_0,{\bf{u}}^*),
\end{aligned}
\end{equation}
where $\beta_s(\theta,S_0)=\frac{1}{S_0}\int_{\max\{0,s-S_0\}}^s d\theta(t), \ \forall s\geq 0$, and $\varsigma(\theta, S_0)$ is the evaluation with $s\mapsto \beta_s(\theta,S_0)$ its density function.

Next, we show that
\begin{align}\label{eq:payoff-dif-S_0}|\gamma_\theta(y_0,{\bf{u}}^*)-\gamma_{\varsigma(\theta,S_0)}(y_0,{\bf{u}}^*)| \leq \sup_{Q\in \mathcal{B}(\RR_+)}|\theta(Q)-\varsigma(\theta,S_0)(Q)|\leq 2TV_{S_0}(\theta).\end{align}
Indeed, the first inequality follows from Hahn's decomposition theorem applied to the sign measure "$\theta-\varsigma(\theta,S_0)$" \big(cf. Lemma 3.7 in \cite{Li_2016}\big). Let $Q$ be any Borel set on $\RR_+$. We write $\beta_s(\theta,S_0)=\frac{1}{S_0}\int_{s-S_0}^s d\theta(t)$ for all $s\geq 0$ by considering $\theta$ as a probability measure over $[-S_0,0)\cup \RR_+$ null on $[-S_0,0)$. We have
\begin{eqnarray*}
\varsigma(\theta,S_0)(Q)=\int_{s\in Q}\beta_s(\theta,S_0) d s&=&\frac{1}{S_0}\int_{s\in Q}\left( \int_{t\in[s-S_0,s]}d\theta(t)\right)d s\\
(''\text{Fubini's theorem}'') &=&\int_{t\in Q-S_0}\left(\frac{1}{S_0}\int_{s\in[t,t+S_0]} d s \right)d\theta(t)\\
&=& \theta(Q-S_0).
\end{eqnarray*}
Thus we have $|\theta(Q)-\varsigma(\theta,S_0)(Q)|=|\theta(Q)-\theta(Q-S_0)|\leq \theta([0,S_0))+TV_{S_0}(\theta)\leq 2TV_{S_0}(\theta)$. This proves Equation (\ref{eq:payoff-dif-S_0}) by taking the supremum over $Q\in \mathcal{B}(\RR_+)$.

Finally, we substitute Equation (\ref{eq:uniform-2}) into Equation (\ref{eq:payoff-dif-S_0}), to obtain:
$$\gamma_{\theta}(y_0,{\bf{u}}^*)\leq  V^*(y_0)+ 2TV_{S_0}(\theta)\leq V^*(y_0)+3\varepsilon,$$
for all $\theta\in\De(\RR^+)$ with $\sup_{0\leq s\leq S_0}TV_s(\theta)\leq \varepsilon$.

To conclude, we have obtained that: $\forall \varepsilon, \exists S_0>0, \exists {\bf{u}}^*\in\tilde{U}$,
$$\forall \theta\in \De(\RR^+),\  \Big(\sup_{0\leq s\leq S_0}TV_s(\theta)\leq \varepsilon\Longrightarrow \gamma_\theta(y_0,{\bf{u}}^*)\leq V^*(y_0)+3\varepsilon,\ \forall y_0\in Y\Big).$$
As $\varepsilon>0$ is arbitrary, this proves Theorem \ref{thm:guv} by taking "$\eta=\varepsilon$" and "$S=S_0$". \\

\section*{Acknowledgments}
\  The main results in this article forms Chapter 4 of my PhD thesis submitted to Universit\'e Pierre et Marie Curie - Paris 6, 2015 September. I wish to thank my supervisor Sylvain Sorin for numerous helpful comments. Discussions with Marc Quincampoix, Fabien Gensbittel and Marco Mazzola are also acknowledged. Thanks go to Cheng Wan for her careful proofreading on part of this manuscript. Part of the work was done when the author was an ATER fellow at Universit\'e Paris-1 (UFR 06) during the academic year 2014-2015, and at Universit\'e Cergy-Pontoise (THEMA $\&$ UFR \'economie et gestion) during the academic year 2015-2016.


\begin{thebibliography}{1}
\bibitem{Alvarez_2010}  {\sc  O. Alvarez and M. Bardi}, {\em Ergodicity, stabilization, and singular pertubations for Bellman-Isaacs equations}, Mem. Amer. Math. Soc., (204)2010,  vi+77 pp.

\bibitem{Arisawa_97} {\sc  M. Arisawa}, {\em Ergodic problem for the Hamilton-Jacobi-Belmann equations}, Ann. Henri Poincar\'e, Analyse Nonlin\'eaire, (14)1997, pp. 415--438.

\bibitem{Arisawa_98} {\sc  M. Arisawa and P. L. Lions}, {\em On ergodic stochastic control}, Comm. Partial Differential Equations, (23)1998, pp. 2187--2217.
 
 \bibitem{Aumann_64} {\sc R. Aumann}, {\em Mixed and behavior strategies in infinite extensive games},  in Advances in Game Theory, Annals of Mathematics Studies 52, M. Dresher and L.S. Shapley eds., 1964, pp. 627--650.

\bibitem{Bensoussan_88} {\sc A. Bensoussan}, {\em Perturbation Methods in Optimal Control}, Wiley/Gauthiers-Villas, Chichester, 1988.

\bibitem{Buckdahn_2014} {\sc R. Buckdahn, D. Goreac and M. Quincampoix}, {\em Existence of asymptotic values for nonexpansive stochastic control systems}, Applied Mathematics and Optimization, (70)2014, pp. 1--28.


\bibitem {Cannarsa_2015} {\sc P. Cannarsa P. and M. Quincampoix}, {\em Vanishing discount limit and nonexpansive optimal control and differential games}, SIAM Journal on Control and Optimization, 53(2007), pp. 1789--1814.
 
\bibitem{Cardaliaguet_2007} {\sc P. Cardaliaguet}, {\em Differential games with asymmetric information}, SIAM Journal on Control and Optimization, 46(2007), pp. 818--838.

\bibitem{Carlson_91} {\sc D. Carlson, A. Haurie and A. Leizarowitz}, {\em  Infinite Horizon Optimal Control: Deterministic and Stochastic Systems}, 2nd ed., Springer, 1991. 

\bibitem {Filar_92} {\sc J. Filar and O.J. Vrieze}, {\em Weighted reward criteria in competitive Markov decision processes}, Zeitschrift fur Operations Research, 36(1992), pp. 343--358. 

\bibitem{Gaitsgory_86} {\sc  V. Gaitsgory}, {\em On the use of the averaging method in control problems}, (Russian) \emph{Differentsialnye Uravneniya}, 22(1986), pp. 1876--1886.

\bibitem {Goreac_2015} {\sc D. Goreac}, {\em  Asymptotic control for a class of piecewise deterministic markov processes associated to temperate viruses}, SIAM J. Control Optimization, 53(2015), pp. 1860--1891.

\bibitem{Goreac_2015b} {\sc  D. Goreac}, {\em A note on general Tauberian-type results for controlled stochastic dynamics},   Electronic Communications in Probability, 20(2015), pp. 1--12.

\bibitem {Haurie_2002} {\sc A. Haurie}, {\em Turnpikes in multidiscount rate environments and GCC policy evaluation}, in Optimal Control and Differential Games, G.Zaccour (ed.), Kluwer, (2002), pp. 39--52. 

\bibitem{Khlopin_2015} {\sc  D. Khlopin}, {\em On uniform Tauberian theorems for dynamic games}, arXiv{:}1412.7331, 2015.


\bibitem{Li_2016} {\sc X. Li, M. Quincampoix and J. Renault}, {\em Limit value for optimal control with general means}, Discrete and Continuous Dynamical Systems--Series A, (36)2016, pp. 2113--2132.

\bibitem{Oliu_2013} {\sc M. Oliu-Barton and G. Vigeral}, {\em A uniform {T}auberian theorem in optimal control}, in Advances in Dynamic Games, Annals of the International Society of Dynamic Games, P. Cardaliaguet and R. Cressman eds., Birkhauser, 12(2013), pp. 199--215.

\bibitem{Quincampoix_2011} {\sc  M. Quincampoix and J. Renault}, {\em On the existence of a limit value in some nonexpansive optimal control problems}, SIAM Journal on Control and Optimization, 49(2011), pp. 2118--2132.

\bibitem{Quincampoix_2015} {\sc M. Quincampoix and S. Hayk}, {\em Averaging problem for weakly coupled nonexpansive control systems}, Nonlinear Analysis: Theory, Methods \& Applications, 113(2015), pp. 147--158.  

\bibitem{Renault_2013} {\sc  J. Renault and X. Venel}, {\em A distance on some probability spaces, with applications to Markov decision processes and repeated games}, hal{:}00674998, 2013.

\bibitem{Sorin_2002} {\sc S. Sorin}, {\em A First Course on Zero-Sum Repeated Games}, Springer, 2002.

\bibitem{Villani_2009} {\sc C. Villani}, {\em Optimal Transportation: Old and New}, Springer, 2009.

\bibitem{Wagner_77} {\sc  D.H. Wagner}, {\em Survey of measurable selection theorems}, SIAM Journal on Control and Optimization, 15(1977), pp. 859--903.
\end{thebibliography}
\end{document}